\documentclass{article}
\usepackage{amsbsy,amssymb,amscd,amsfonts,latexsym,amstext,delarray,
amsmath, diagrams}  \headsep=15pt
\usepackage[all,knot]{xy}
\def\cX{{\cal X} }

\def\bs{ \backslash}

\def\G{{ \Gamma }}

\def\cE{ {\cal E}}

\def\bq{\begin{quote}}
\def\eq{\end{quote}}

\def\P{ {\bf P}}

\def\be{\begin{equation}}
\def\ee{\end{equation}}

\def\bd{\begin{diagram} }
\def\ed{ \end{diagram}}

\def\bs{\begin{slide}}
\def\es{\end{slide} }

\def\bd{ \begin{diagram}}
\def\ed{\end{diagram} }

\author{Minhyong Kim}
\begin{document}
\begin{flushleft}
{\bf Non-abelian fundamental groups in arithmetic geometry}
\medskip

{\bf Minhyong Kim}
\end{flushleft}

\medskip

The entry on  mathematics in the Wikipedia (as of evening, 28 September, 2009)
characterizes it as  `the science and study of quantity, structure, space, and change.'
Arithmetic geometry reflects with great intensity on the first three,
in that a remarkable {\em unity} of space and quantity is effected through
the medium of highly abstract, but precise and robust structures.
On the other hand, it is indeed an area of mathematics in which
{\em time}  plays almost no role, separating it from the parts
of mathematics traditionally connected to the palpable phenomena.

That quantities make up the concrete language of
space should surprise no one, inasmuch as distances, heights, and other commonplace dimensions
are described in numbers even by the innocent bystander.
And then,  subtle {\em relations and constraints} on a collection of quantities might hint at
 an intricacy of geometric constitution, as in the collection of latitudes and
longitudes for positions on our globe. Closer to the textbooks,
relations between `abstract quantities' indicated by funny variables like `$x$' and `$y$',
as when we write
$$x^2+y^2=1,$$
can reflect geometry
of a very harmonious nature. Perhaps some will even remember how to trace the pairs $(x,y)$ in the plane
that are so related, to recover the physical form that realizes the abstract one.

A  deep phenomenon, exceedingly powerful and far-reaching
in spite of its philosophical garb,
is a beautiful {\em duality} between space and quantity, or geometry and algebra, whose
basic properties were
articulated most clearly by Alexander Grothendieck in the 1960's, and came to
be known as the theory of {\em schemes}. Almost any algebraic problem can be fashioned into
a geometric picture, with far greater fluency, in fact, than the vice versa. Even a seasoned practitioner
 is often enough mystified
that such a geometric framework underlies humble equations in two variables
as well as, on occasion, the fabric of the physical universe itself\footnote{It should be furthermore noted
that the age of computers has brought forth the importance of
arithmetic geometric structure in many problems of  packaging and processing
information in an efficient fashion. In more fanciful terms, one might say that the {\em artificial universe}
has spectacularly incorporated arithmetic geometry.}.

The fact that one can visualize a solution to an equation, say
$$(5/13)^2+(12/13)^2=1,$$
as a {\em point}, is already indicative of a very general procedure for extracting space out of
algebra. This particular equation
has an infinity of solutions in rational numbers that we can collect by
arranging  pairs of whole numbers $(m,n)$ into the form
$$(\frac{m^2-n^2}{m^2+n^2}, \frac{2mn}{m^2+n^2}).$$
On the other hand, the equations
$$3x^4+5y^4=6$$
and
$$y^2=x^5 - 14x^4 + 65x^3 - 112x^2 + 60x$$
both have only finitely many rational solutions, while
$$y^2=2374618x^5 - 44158534x^4 + 81193214x^3 - 39409298x^2
$$
and
$$x^3+y^3=1729$$
have again  infinite solution sets. For the equation
$y^2=x^5 - 14x^4 + 65x^3 - 112x^2 + 60x$, in fact, the only  solutions are
$$(0,0), (1,0), (2,0), (5,0), (6,0), (3,6), (3,-6), (10,120), (10,-120).$$
Here is a rather difficult solution to the cubic equation:
$$(\frac{-5150812031}{107557668})^3+(\frac{5177701439}{107557668})^3=1729,$$
and one can systematically find many more, even though the method for generating them is a bit more complex than for
the circle.

This pattern of difference, and the question of {\em which}  equations in two variables
will have finitely many solutions, is almost well-understood, thanks to a celebrated theorem
of Gerd Faltings.
What lies beyond our present scope is a substantive grip on the
{\em ultimate reason} generating these differences. Plausible attempts at explanation
are as numerous in mathematics as in the efforts of normal science to grapple with
 core natural phenomena. In the context of our current programme, the goal is to analyze
the pattern completely in terms of the {\em fundamental group}
of the equation, a highly structured measure of its complexity.
Here, critical use is made of the space  associated to an equation, in that the fundamental group
 weighs the totality of  paths that one might attempt to traverse through it.

An important (indeed `fundamental') distinction
is between abelian and non-abelian fundamental groups, where the latter
require of us  vastly greater numbers of  distinct paths to get from one point to another
than does the former. That the complexity of  non-abelian structures, and the eventual generality of {\em higher-dimensional algebra}, might be refined into  constructive machinery
for solving every possible equation was the vision of Grothendieck emanating from the 1980's, a disquieting decade
that ended  with his  disappearance into the Pyrenees. Yet, this strange proposal  motivates our present
gathering at the Newton Institute.
\medskip

\begin{flushleft}
{\bf Appendix}

\medskip

We wish here to convey a feeling for non-abelian fundamental groups, without being precise
about any of the words that go into that phrase. In mathematics, it is common to refer to any object considered from a geometric viewpoint  as a {\em space}, and we shall do so here. There are  two spaces we will consider.
\smallskip

The first one, which I will refer to as $\cE$ can be constructed  from a sheet of rubber
as follows.
$$
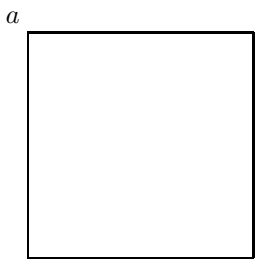
\xy
(-2, 32)*{a};
(0,0)*{}; (30,0)*{} **\dir{-};
(30,0)*{}; (30,30)*{} **\dir{-};
(30,30)*{}; (0,30)*{} **\dir{-};
(0,30)*{}; (0,0)*{} **\dir{-}
\endxy
$$
Step one:

Glue the top edge to the bottom edge. I hope you can see that what you have is a cylinder.
You will see also that the two side edges have now curled up into circles.
\smallskip

Step two:

Bring  those two circles together and glue them. That's it! We're done with this space.
\smallskip

The space $\cE$ you have is called a {\em torus} in mathematical parlance. In reality,
it would look like the inner tube of a tire. Since my skills
for drawing such objects using this editing program is extremely limited,
I will continue to describe the space using the square picture above. You just need to remember
that the bottom and top edges are glued together, as are the two side edges.
Of course you are used to this kind of thing from looking at a world map.
Even though we lay it out flat, you understand that if you go into the right edge, then you come out from the left in real life. Just for practice, ask yourself what kind of path is described by the line $L$ in the following picture:
$$
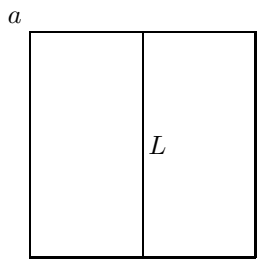
\xy
(-2, 32)*{a};
(17,15)*{L};
(15,30)*{}; (15,0)*{} **\dir{-};
(0,0)*{}; (30,0)*{} **\dir{-};
(30,0)*{}; (30,30)*{} **\dir{-};
(30,30)*{}; (0,30)*{} **\dir{-};
(0,30)*{}; (0,0)*{} **\dir{-}
\endxy
$$
I'm sure you've guessed that it's actually a circle.
\smallskip

For something slightly trickier, try
$$
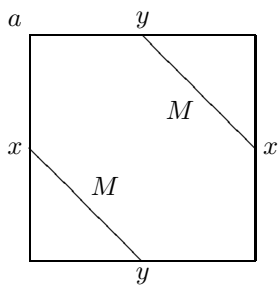
\xy
(-2,15)*{x};
(32,15)*{x};
(15, -2)*{y};
(15,32)*{y};
(-2, 32)*{a};
(10,10)*{M};
(20,20)*{M};
(0,15)*{}; (15,0)*{} **\dir{-};
(15,30)*{}; (30,15)*{} **\dir{-};
(0,0)*{}; (30,0)*{} **\dir{-};
(30,0)*{}; (30,30)*{} **\dir{-};
(30,30)*{}; (0,30)*{} **\dir{-};
(0,30)*{}; (0,0)*{} **\dir{-}
\endxy
$$
Notice that it looks like there are two segments, both of which I've labeled
$M$. But on the torus, you'll see that you've ended up again with a single circle,
because the two points labeled $y$ are actually joined, as are the points labeled $x$.
\smallskip

On the other hand, if we look at something like
$$
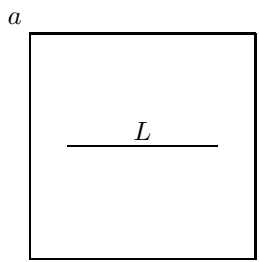
\xy
(-2, 32)*{a};
(15,17)*{L};
(5,15)*{}; (25,15)*{} **\dir{-};
(0,0)*{}; (30,0)*{} **\dir{-};
(30,0)*{}; (30,30)*{} **\dir{-};
(30,30)*{}; (0,30)*{} **\dir{-};
(0,30)*{}; (0,0)*{} **\dir{-}
\endxy
$$
It will clearly remain a segment even on the torus. The ends just don't join up even after we've
glued the edges as instructed.

$$
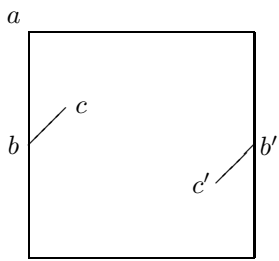
\xy
(-2, 32)*{a};
(-2,15)*{b};
(32,15)*{b'};
(7,20)*{c};
(23,10)*{c'};
(0,15)*{}; (5,20)*{} **\dir{-};
(25,10)*{}; (30,15)*{} **\dir{-};
(0,0)*{}; (30,0)*{} **\dir{-};
(30,0)*{}; (30,30)*{} **\dir{-};
(30,30)*{}; (0,30)*{} **\dir{-};
(0,30)*{}; (0,0)*{} **\dir{-}
\endxy
$$
For the line above, you will see that the points $b$ and $b'$ are the same on
the torus, but $c$ and $c'$ remain distinct. So the result is a single line segment.
\smallskip

I spoke above about the number of ways to get from one point to another in some space.
The counting in this context has to be done according to specific conventions, and we will move towards
a discussion of them.
\smallskip

First of all, consider the following path $P$, emanating from the point labeled $a$. It's probably best to think about walking right along the edge, even though I've drawn it just below for ease of viewing.
$$
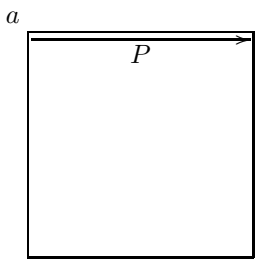
\xy
(15,27)*{P};
{\ar (0.5,29)*{}; (29.5,29)*{}};
(-2, 32)*{a};
(0,0)*{}; (30,0)*{} **\dir{-};
(30,0)*{}; (30,30)*{} **\dir{-};
(30,30)*{}; (0,30)*{} **\dir{-};
(0,30)*{}; (0,0)*{} **\dir{-}
\endxy
$$
Where does $P$ end? Remember that our square is just a map of a torus.
In particular, since the two side edges are glued together, the end point of $P$ is just...
$a$! Like the line $L$, the path $P$ will simply go around a circle on the torus $E$.
It's good to keep in mind that the square map and the torus itself are pretty easy to visualize separately.
The tricky part is describing the geometry on the torus in terms of the square map. Why do we
do this? It's essentially because explaining  on flat paper is considerably easier than
carrying around balloons or tires. Even if I did, it would be quite cumbersome to draw diagrams
on them and keep track of the different portions as we rotate them around in the course of a discussion.
If you become a great engineer in the future, you should invent
a 3D box in which these operations can be carried out as easily as we now write in our notebook. It is
plausible that the kind of `methodology of geometric description' we are in the middle of right now will
be an important component of the technology, especially since there needs to be an efficient
interface with an embedded computer, which will need to have everything spelled out in numbers.
Once the technology for discussing 3D is sufficiently developed, we may be able to work with
4D and higher with far greater ease than we can now.
$$
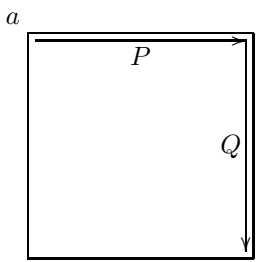
\xy
(15,27)*{P};
{\ar (1,29)*{}; (29,29)*{}};
(27,15)*{Q};
{\ar (29,29)*{}; (29,1)*{}};
(-2, 32)*{a};
(0,0)*{}; (30,0)*{} **\dir{-};
(30,0)*{}; (30,30)*{} **\dir{-};
(30,30)*{}; (0,30)*{} **\dir{-};
(0,30)*{}; (0,0)*{} **\dir{-}
\endxy
$$
If we look at the path $Q$, I hope it's obvious by now that's it's also a circle, as in the case of the line $L$, because
the top and bottom edges are also glued.
$$
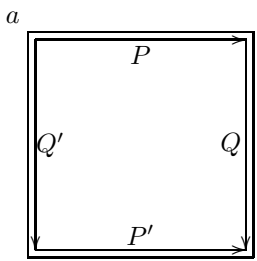
\xy
(15,27)*{P};
{\ar (1,29)*{}; (29,29)*{}};
(27,15)*{Q};
{\ar (29,29)*{}; (29,1)*{}};
(3,15)*{Q'};
{\ar (1,29)*{}; (1,1)*{}};
(15,3)*{P'};
{\ar (1,1)*{}; (29,1)*{}};
(-2, 32)*{a};
(0,0)*{}; (30,0)*{} **\dir{-};
(30,0)*{}; (30,30)*{} **\dir{-};
(30,30)*{}; (0,30)*{} **\dir{-};
(0,30)*{}; (0,0)*{} **\dir{-}
\endxy
$$
Now I've indicated paths $P$, $P'$, $Q$, $Q'$. But the point is that on the torus, there are only two.
Therefore, I've labeled all of them with the same letter in the diagram below.
One other point to note is that I've labeled all four vertices of our square with a single letter
$a$, in honor of the fact that all four points represent the {\em same} point on the torus. Try to make yourself certain of this.
$$
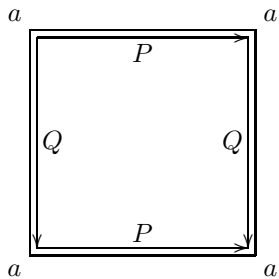
\xy
(-2,-2)*{a};
(32,-2)*{a};
(32,32)*{a};
(15,27)*{P};
{\ar (1,29)*{}; (29,29)*{}};
(27,15)*{Q};
{\ar (29,29)*{}; (29,1)*{}};
(3,15)*{Q};
{\ar (1,29)*{}; (1,1)*{}};
(15,3)*{P};
{\ar (1,1)*{}; (29,1)*{}};
(-2, 32)*{a};
(0,0)*{}; (30,0)*{} **\dir{-};
(30,0)*{}; (30,30)*{} **\dir{-};
(30,30)*{}; (0,30)*{} **\dir{-};
(0,30)*{}; (0,0)*{} **\dir{-}
\endxy$$
To summarize, what appears on our map as four distinct lines are actually just two circles $P$ and $Q$
that meet at the point $a$.
\smallskip

Now to the question of counting paths. We are in fact interested in counting paths between any two fixed points.
But it turns out to be simpler to count paths that end and begin at the same point. Perhaps I will explain sometime
the relation between the two counts. In some sense they are exactly the same. For the moment, we will forget about
that issue and simply count the number of paths that both begin {\em and} end at the point $a$.
These are commonly  referred to as {\em loops based at $a$}.
\smallskip

An essential point that we will accept casually is captured by the following rule:
\bq
We regard two paths $\G_1$ and $\G_2$ as the same if one can be deformed into the
other without moving the endpoints.
\eq
The fancy terminology is that $\G$ and $\G'$ are {\em homotopic}.
For example, the two paths below are homotopic
$$
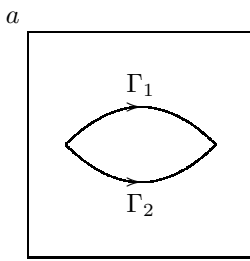
\xy
(-2, 32)*{a};
(15,23)*{\G_1};
(15,7)*{\G_2};
(5,15)*{}; (25,15)*{} **\crv{(15,25)}?(.5)*\dir{>};
(5,15)*{}; (25,15)*{} **\crv{(15,5)}?(.5)*\dir{>};
(0,0)*{}; (30,0)*{} **\dir{-};
(30,0)*{}; (30,30)*{} **\dir{-};
(30,30)*{}; (0,30)*{} **\dir{-};
(0,30)*{}; (0,0)*{} **\dir{-}
\endxy
$$
as one can see by deforming one to the other using a family of intermediate paths:
$$
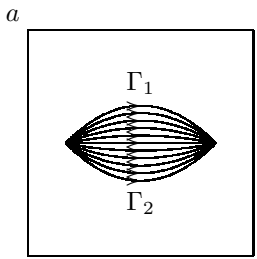
\xy
(-2, 32)*{a};
(15,23)*{\G_1};
(15,7)*{\G_2};
(5,15)*{}; (25,15)*{} **\crv{(15,25)}?(.5)*\dir{>};
(5,15)*{}; (25,15)*{} **\crv{(15,23)}?(.5)*\dir{>};
(5,15)*{}; (25,15)*{} **\crv{(15,21)}?(.5)*\dir{>};
(5,15)*{}; (25,15)*{} **\crv{(15,19)}?(.5)*\dir{>};
(5,15)*{}; (25,15)*{} **\crv{(15,17)}?(.5)*\dir{>};
(5,15)*{}; (25,15)*{} **\crv{(15,15)}?(.5)*\dir{>};
(5,15)*{}; (25,15)*{} **\crv{(15,13)}?(.5)*\dir{>};
(5,15)*{}; (25,15)*{} **\crv{(15,11)}?(.5)*\dir{>};
(5,15)*{}; (25,15)*{} **\crv{(15,9)}?(.5)*\dir{>};
(5,15)*{}; (25,15)*{} **\crv{(15,7)}?(.5)*\dir{>};
(5,15)*{}; (25,15)*{} **\crv{(15,5)}?(.5)*\dir{>};
(0,0)*{}; (30,0)*{} **\dir{-};
(30,0)*{}; (30,30)*{} **\dir{-};
(30,30)*{}; (0,30)*{} **\dir{-};
(0,30)*{}; (0,0)*{} **\dir{-}
\endxy
$$
The family of intermediate paths themselves are referred to as the {\em homotopy},
in this case between $\G_1$ and $\G_2$.
You should be able to convince yourself that both $\G_1$ and $\G_2$ are homotopic to
$\G_3$ below.
$$
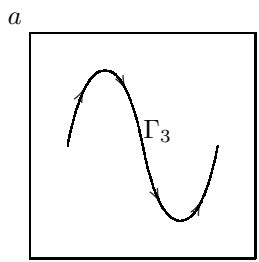
\xy
(-2, 32)*{a};
(17,17)*{\G_3};
(5,15)*{}; (10,25)*{}**\crv{(7,25)}?(.5)*\dir{>};
(10,25)*{}; (15,15)*{} **\crv{(13,25)}?(.5)*\dir{>};
(15,15)*{}; (20,5)*{} **\crv{(17,5)}?(.5)*\dir{>};
(20,5)*{}; (25,15)*{} **\crv{(23,5)}?(.5)*\dir{>};
(0,0)*{}; (30,0)*{} **\dir{-};
(30,0)*{}; (30,30)*{} **\dir{-};
(30,30)*{}; (0,30)*{} **\dir{-};
(0,30)*{}; (0,0)*{} **\dir{-}
\endxy
$$
All three are also homotopic to the piecewise linear path $\G_4$:
$$
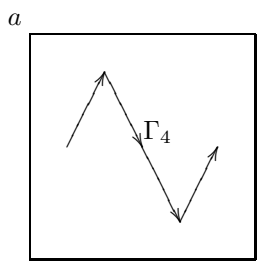
\xy
(-2, 32)*{a};
(17,17)*{\G_4};
{\ar (5,15)*{}; (10,25)*{}};
{\ar (10,25)*{}; (15,15)*{}};
{\ar (15,15)*{}; (20,5)*{}};
{\ar (20,5)*{}; (25,15)*{}};
(0,0)*{}; (30,0)*{} **\dir{-};
(30,0)*{}; (30,30)*{} **\dir{-};
(30,30)*{}; (0,30)*{} **\dir{-};
(0,30)*{}; (0,0)*{} **\dir{-}
\endxy
$$
We will see some more interesting examples of homotopies later on, but it is important to know that
there are many collections of paths that are {\em not} homotopic to each other. For example, it is a fact that
the two paths
$P$ and $Q$ are not homotopic.
$$
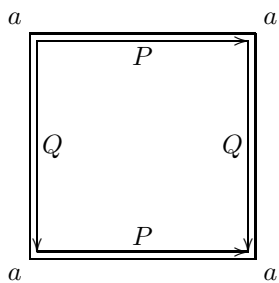
\xy
(15,27)*{P};
(-2,-2)*{a};
(32,-2)*{a};
(32,32)*{a};
{\ar (1,29)*{}; (29,29)*{}};
(27,15)*{Q};
{\ar (29,29)*{}; (29,1)*{}};
(3,15)*{Q};
{\ar (1,29)*{}; (1,1)*{}};
(15,3)*{P};
{\ar (1,1)*{}; (29,1)*{}};
(-2, 32)*{a};
(0,0)*{}; (30,0)*{} **\dir{-};
(30,0)*{}; (30,30)*{} **\dir{-};
(30,30)*{}; (0,30)*{} **\dir{-};
(0,30)*{}; (0,0)*{} **\dir{-}
\endxy
$$
You might, for example, want to deform one into the other
using a family like this:
$$
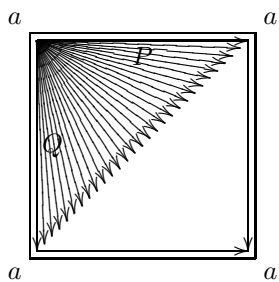
\xy
(15,27)*{P};
(-2,-2)*{a};
(32,-2)*{a};
(32,32)*{a};
{\ar (1,29)*{};(2,2)*{}};
{\ar (1,29)*{};(3,3)*{}};
{\ar (1,29)*{};(4,4)*{}};
{\ar (1,29)*{};(5,5)*{}};
{\ar (1,29)*{};(6,6)*{}};
{\ar (1,29)*{};(7,7)*{}};
{\ar (1,29)*{};(8,8)*{}};
{\ar (1,29)*{}; (9,9)*{}};
{\ar (1,29)*{}; (10,10)*{}};
{\ar (1,29)*{}; (11,11)*{}};
{\ar (1,29)*{}; (12,12)*{}};
{\ar (1,29)*{}; (13,13)*{}};
{\ar (1,29)*{}; (14,14)*{}};
{\ar (1,29)*{}; (15,15)*{}};
{\ar (1,29)*{}; (16,16.)*{}};
{\ar (1,29)*{}; (17,17)*{}};
{\ar (1,29)*{}; (18,18)*{}};
{\ar (1,29)*{}; (19,19)*{}};
{\ar (1,29)*{}; (20,20)*{}};
{\ar (1,29)*{}; (21,21)*{}};
{\ar (1,29)*{}; (22,22)*{}};
{\ar (1,29)*{}; (23,23)*{}};
{\ar (1,29)*{}; (24,24)*{}};
{\ar (1,29)*{}; (25,25)*{}};
{\ar (1,29)*{}; (26,26.)*{}};
{\ar (1,29)*{}; (27,27)*{}};
{\ar (1,29)*{}; (28,28)*{}};
{\ar (1,29)*{};(29,29)*{}};
{\ar (29,29)*{}; (29,1)*{}};
(3,15)*{Q};
{\ar (1,29)*{}; (1,1)*{}};
{\ar (1,1)*{}; (29,1)*{}};
(-2, 32)*{a};
(0,0)*{}; (30,0)*{} **\dir{-};
(30,0)*{}; (30,30)*{} **\dir{-};
(30,30)*{}; (0,30)*{} **\dir{-};
(0,30)*{}; (0,0)*{} **\dir{-}
\endxy
$$
but this would not a homotopy. The rule, remember, is that two paths that can be deformed into one another
are homotopic only if the deformation leaves the beginning and end points fixed. It is actually not so easy to
prove beyond doubt that $P$ and $Q$ are not homotopic. The clear proof is usually
regarded as advanced undergraduate level mathematics, and we will not go further into it.
Instead, we will look at a legitimate homotopy from $P$ to some other paths.
We can start by just pulling it down a little bit without moving the endpoints.
$$
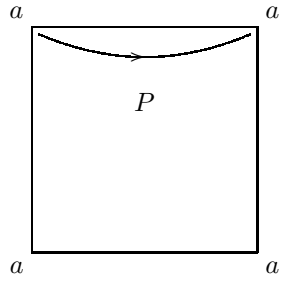
\xy
(-2,-2)*{a};
(32,-2)*{a};
(32,32)*{a};
(15,20)*{P};
(1,29)*{}; (29,29)*{} **\crv{(15,23)}?(.5)*\dir{>};
(-2, 32)*{a};
(0,0)*{}; (30,0)*{} **\dir{-};
(30,0)*{}; (30,30)*{} **\dir{-};
(30,30)*{}; (0,30)*{} **\dir{-};
(0,30)*{}; (0,0)*{} **\dir{-}
\endxy
$$
Then we keep going.
$$
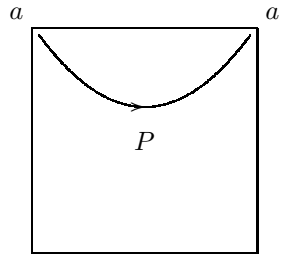
\xy
(-2,-2)*{a};
(32,-2)*{a};
(32,32)*{a};
(15,15)*{P};
(1,29)*{}; (29,29)*{} **\crv{(15,10)}?(.5)*\dir{>};
(-2, 32)*{a};
(0,0)*{}; (30,0)*{} **\dir{-};
(30,0)*{}; (30,30)*{} **\dir{-};
(30,30)*{}; (0,30)*{} **\dir{-};
(0,30)*{}; (0,0)*{} **\dir{-}
\endxy
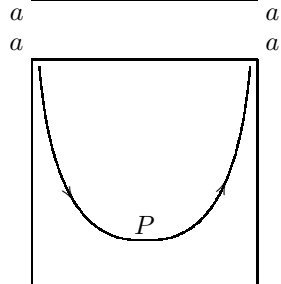
$$
$$
\xy
(-2,-2)*{a};
(32,-2)*{a};
(32,32)*{a};
(15,8)*{P};
(1,29)*{}; (15,6)*{} **\crv{(3,5)}?(.5)*\dir{>};
(15,6)*{}; (29,29)*{} **\crv{(27,5)}?(.6)*\dir{>};
(-2, 32)*{a};
(0,0)*{}; (30,0)*{} **\dir{-};
(30,0)*{}; (30,30)*{} **\dir{-};
(30,30)*{}; (0,30)*{} **\dir{-};
(0,30)*{}; (0,0)*{} **\dir{-}
\endxy
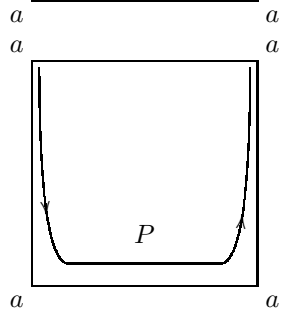
$$
$$
\xy
(-2,-2)*{a};
(32,-2)*{a};
(32,32)*{a};
(15,7)*{P};
(1,29)*{}; (5,3)*{} **\crv{(1,3)}?(.5)*\dir{>};
(5,3)*{}; (25,3)*{} **\dir{-};
(25,3)*{}; (29,29)*{} **\crv{(29,3)}?(.5)*\dir{>};
(-2, 32)*{a};
(0,0)*{}; (30,0)*{} **\dir{-};
(30,0)*{}; (30,30)*{} **\dir{-};
(30,30)*{}; (0,30)*{} **\dir{-};
(0,30)*{}; (0,0)*{} **\dir{-}
\endxy
$$
It shouldn't be hard to convince yourself that we end up thereby with a homotopy from $P$ to the
path that goes along $Q$, then $P$, and then {\em backwards} along $Q$. We show it below, with
the backward route along $Q$ labeled suggestively by $Q^{-1}$.
$$
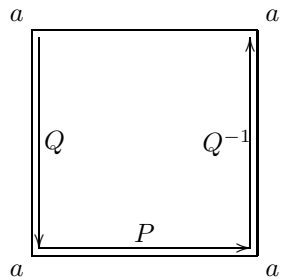
\xy
(-2,-2)*{a};
(32,-2)*{a};
(32,32)*{a};
(26,15)*{Q^{-1}};
{\ar (29,1)*{}; (29,29)*{}};
(3,15)*{Q};
{\ar (1,29)*{}; (1,1)*{}};
(15,3)*{P};
{\ar (1,1)*{}; (29,1)*{}};
(-2, 32)*{a};
(0,0)*{}; (30,0)*{} **\dir{-};
(30,0)*{}; (30,30)*{} **\dir{-};
(30,30)*{}; (0,30)*{} **\dir{-};
(0,30)*{}; (0,0)*{} **\dir{-}
\endxy
$$
Since we have agreed to regard homotopic paths as the {\em same}, we can express the
previous conclusion as a mathematical equation:
$$P=QPQ^{-1},$$
where the path on the right refers to the concatenation of the three that make it up: $Q$, then
$P$, then $Q$ backwards. Here, we take a moment to note that
this path is made of three others which all begin and end at $a$. That is, we start at $a$
and come back three times, but consider the whole thing as one path. For another example, $PP$,
going along the same path $P$ twice, is not necessarily\footnote{
We can't immediately rule out the possibility that they might be homotopic, although they end up not to be in this case.} to be regarded as the same as $P$.
Of course, in real life as well, running the same circuit twice certainly feels different from
doing it once. But here, we are regarding some paths as the same, i.e., when they are
homotopic. So I guess you do need to trust me somewhat that all these different
conventions are somehow consistent and interesting.

\smallskip
We will now display a very important homotopy, first in a few stages:
$$
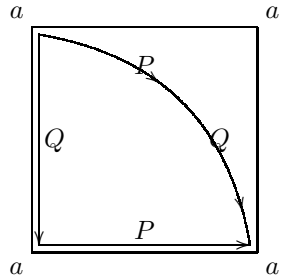
\xy
(-2,-2)*{a};
(32,-2)*{a};
(32,32)*{a};
(15,25)*{P};
(25,15)*{Q};
(1,29)*{}; (29,1)*{} **\crv{(25,25)}?(.4)*\dir{>}?(.9)*\dir{>};
(3,15)*{Q};
{\ar (1,29)*{}; (1,1)*{}};
(15,3)*{P};
{\ar (1,1)*{}; (29,1)*{}};
(-2, 32)*{a};
(0,0)*{}; (30,0)*{} **\dir{-};
(30,0)*{}; (30,30)*{} **\dir{-};
(30,30)*{}; (0,30)*{} **\dir{-};
(0,30)*{}; (0,0)*{} **\dir{-}
\endxy
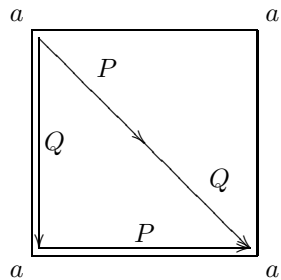
$$
$$
\xy
(-2,-2)*{a};
(32,-2)*{a};
(32,32)*{a};
(10,25)*{P};
(25,10)*{Q};
{\ar (1,29)*{};(15,15)*{}};
{\ar (15,15)*{};(29,1)*{}};
(3,15)*{Q};
{\ar (1,29)*{}; (1,1)*{}};
(15,3)*{P};
{\ar (1,1)*{}; (29,1)*{}};
(-2, 32)*{a};
(0,0)*{}; (30,0)*{} **\dir{-};
(30,0)*{}; (30,30)*{} **\dir{-};
(30,30)*{}; (0,30)*{} **\dir{-};
(0,30)*{}; (0,0)*{} **\dir{-}
\endxy
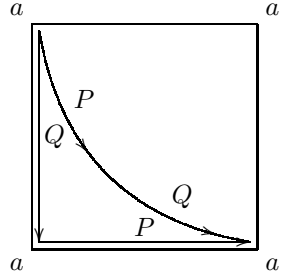
$$
$$
\xy
(-2,-2)*{a};
(32,-2)*{a};
(32,32)*{a};
(7,20)*{P};
(20,7)*{Q};
(1,29)*{}; (29,1)*{} **\crv{(5,5)}?(.4)*\dir{>}?(.9)*\dir{>};
(3,15)*{Q};
{\ar (1,29)*{}; (1,1)*{}};
(15,3)*{P};
{\ar (1,1)*{}; (29,1)*{}};
(-2, 32)*{a};
(0,0)*{}; (30,0)*{} **\dir{-};
(30,0)*{}; (30,30)*{} **\dir{-};
(30,30)*{}; (0,30)*{} **\dir{-};
(0,30)*{}; (0,0)*{} **\dir{-}
\endxy
$$
and then, all at once
$$
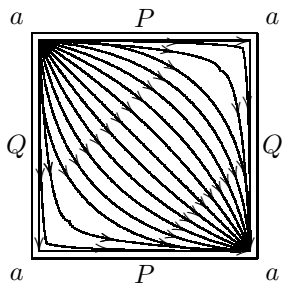
\xy
{\ar (1,29)*{}; (29,29)*{}};
{\ar (29,29)*{}; (29,1)*{}};
(-2,-2)*{a};
(32,-2)*{a};
(32,32)*{a};
(15,32)*{P};
(32,15)*{Q};
(1,29)*{}; (28,28)*{} **\crv{(25,29)}?(.5)*\dir{>};
(28,28)*{};(29,1)*{} **\crv{(29,25)}?(.5)*\dir{>};
(1,29)*{}; (2,2)*{} **\crv{(1,5)}?(.5)*\dir{>};
(2,2)*{};(29,1)*{} **\crv{(5,1)}?(.5)*\dir{>};
(1,29)*{}; (5,5)*{} **\crv{(2,5)}?(.5)*\dir{>};
(5,5)*{};(29,1)*{} **\crv{(5,2)}?(.5)*\dir{>};
(1,29)*{}; (25,25)*{} **\crv{(25,28)}?(.5)*\dir{>};
(25,25)*{};(29,1)*{} **\crv{(28,25)}?(.5)*\dir{>};
(1,29)*{}; (5,5)*{} **\crv{(2,5)}?(.5)*\dir{>};
(5,5)*{};(29,1)*{} **\crv{(5,2)}?(.5)*\dir{>};
(1,29)*{}; (29,1)*{} **\crv{(30,30)}?(.4)*\dir{>}?(.8)*\dir{>};
(1,29)*{}; (29,1)*{} **\crv{(27,27)}?(.4)*\dir{>}?(.8)*\dir{>};
(1,29)*{}; (29,1)*{} **\crv{(24,24)}?(.4)*\dir{>}?(.8)*\dir{>};
(1,29)*{}; (29,1)*{} **\crv{(21,21)}?(.4)*\dir{>}?(.8)*\dir{>};
(1,29)*{}; (29,1)*{} **\crv{(18,18)}?(.4)*\dir{>}?(.8)*\dir{>};
(1,29)*{}; (29,1)*{} **\crv{(15,15)}?(.4)*\dir{>}?(.8)*\dir{>};
(1,29)*{}; (29,1)*{} **\crv{(12,12)}?(.4)*\dir{>}?(.8)*\dir{>};
(1,29)*{}; (29,1)*{} **\crv{(9,9)}?(.4)*\dir{>}?(.8)*\dir{>};
(1,29)*{}; (29,1)*{} **\crv{(6,6)}?(.4)*\dir{>}?(.8)*\dir{>};
(1,29)*{}; (29,1)*{} **\crv{(3,3)}?(.4)*\dir{>}?(.8)*\dir{>};
(1,29)*{}; (29,1)*{} **\crv{(0,0)}?(.4)*\dir{>}?(.8)*\dir{>};
(-2,15)*{Q};
{\ar (1,29)*{}; (1,1)*{}};
(15,-2)*{P};
{\ar (1,1)*{}; (29,1)*{}};
(-2, 32)*{a};
(0,0)*{}; (30,0)*{} **\dir{-};
(30,0)*{}; (30,30)*{} **\dir{-};
(30,30)*{}; (0,30)*{} **\dir{-};
(0,30)*{}; (0,0)*{} **\dir{-}
\endxy
$$
The point is that the path that traverses $P$ and then $Q$, is homotopic to the
path the goes along $Q$ first, and then along $P$.
As an equation, we can write:
$$PQ=QP.$$
This is the key equation that characterizes the fact that the fundamental group (whatever it is) of the space $\cE$
is {\em abelian}. This equation allows us to change orders in the way we go along a path
and end up with a path that is essentially the same. For example, if we consider the path
$$PPQQPQQ,$$
we can use the equation above to exchange the order of the of the second $Q$ and the third $P$, ending up
with
$$PPQPQQQ$$
and then exchange the first $Q$ and the third $P$, leaving us with
$$PPPQQQQ$$
That is,
$$PPQQPQQ=PPPQQQQ.$$
you should convince yourself that this is also the same, for example, as
$$QQQQPPP.$$
\smallskip

It is a fact, again at the level of university mathematics, that
any path in $\cE$ that begins and ends in $a$ is homotopic to a path obtained entirely
by traversing a certain number of $Q$'s and a certain number of $P$'s, in a certain order.
Since we have just shown that we can change the orders all we want, we see that any path
in $\cE$ that starts and ends at $a$ can be written
as
$$PPP\cdots PQQQ\cdots QQ,$$
with some $P$'s at the beginning and then a bunch of $Q$'s.
\smallskip

Suppose we want to list all the paths with endpoints at $a$
having length $4$. Then we might have all $P$'s:
$$PPPP,$$
three $P$'s and one $Q$:
$$PPPQ,$$
two $P$'s and two $Q$'s:
$$PPQQ,$$
one $P$ and three $Q$'s:
$$PQQQ,$$
finally, a path of four $Q$'s:
$$QQQQ.$$
You might protest that we haven't counted
$PQPQ$, for example. But this is the same as $PPQQ$, since
we have already seen that the order of any $P$ and $Q$ can be switched. To spell this out just once more in words,
\bq
To go along $P$, then $Q$, then $P$, then $Q$,
\eq
is the same as
\bq
going along $P$, then the path
$QP$, and then $Q$,
\eq
 which is the same as
 \bq
 $P$, then $PQ$, then $Q$,
 (since $QP$ and $PQ$ are the same).
 \eq
That is,
\bq
$P$, then $P$, then $Q$, then $Q$.
\eq
 I hope you can tolerate that bit of over-explanation.

\smallskip
So how many paths would there be with length 10?
We could start listing them again
$$PPPPPPPPPP$$
$$PPPPPPPPPQ$$
$$...$$
or simply observe that the path will only depend on how many $P$'s and $Q$'s there
are. It eases the problem of counting to make the observation that the number of
$Q$'s is completely determined by the number of $P$'s, since for example, if there
are 3 $P$'s then there must be 7 $Q$'s. With this preliminary remark, we see
that there could be $0,1,2,\ldots, $ up to 10 $P$'s. This gives us 11 possibilities.
(It's a bit confusing that it's not 10 possibilities. The reason is that we start counting at 0.)
\smallskip

We can generalize this, and summarize our many observations with a {\em theorem}:
\bq
On the space $\cE$, for each whole number $n$, there are $n+1$ genuinely distinct
paths of length $n$ that start and end at $a$.
\eq

Of course, I haven't proved this for you. In particular, you need to take on faith that
the no two of the paths listed in the theorem are not somehow homotopic to each other in
a clever, unpredictable way. Nevertheless, I hope our discussion does give you at least some feel
for how this all works.
\smallskip

Now I will explain to you, very briefly, another space that looks quite similar in some ways, but
with a dramatic difference in the count of paths. We will call it $\cX$. The space
$\cX$  starts with the cylinder $\cE$ already discussed, and adds to the construction a
\smallskip

Step 3:

Punch a round hole in the middle.
\smallskip

It looks like this:
$$
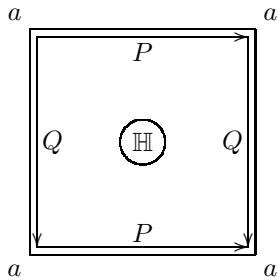
\xy
(-2,-2)*{a};
(32,-2)*{a};
(32,32)*{a};
(15,15)*{ \mathbb{H}};
(15,15)*\xycircle(3,3){-};
(15,27)*{P};
{\ar (1,29)*{}; (29,29)*{}};
(27,15)*{Q};
{\ar (29,29)*{}; (29,1)*{}};
(3,15)*{Q};
{\ar (1,29)*{}; (1,1)*{}};
(15,3)*{P};
{\ar (1,1)*{}; (29,1)*{}};
(-2, 32)*{a};
(0,0)*{}; (30,0)*{} **\dir{-};
(30,0)*{}; (30,30)*{} **\dir{-};
(30,30)*{}; (0,30)*{} **\dir{-};
(0,30)*{}; (0,0)*{} **\dir{-}
\endxy
$$
The edges are all glued as before, but the region I've labeled $\mathbb{H}$
is now missing. The actual object would look like  the inner tube of a tire with a round hole
in the rubber. However, most of the old paths are there, in particular, our friends $P$ and $Q$.
\smallskip

I will not bore you with a long discussion of the properties of $\cX$. Instead,
I will zoom in on the most important one for our purposes. In $\cX$,
$$PQ\neq QP,$$ that is $PQ$ and $QP$ are {\em not} homotopic. But how could that be?
Aren't the paths exactly the same as before? Well they are, but  it is the homotopies that
have changed as a result of moving from $\cE$ to the space $\cX$.
\smallskip

Once again, it is actually not so easy to prove conclusively
that $PQ$ and $QP$ are not homotopic.
But it is instructive to try the previous homotopy from $PQ$ to $QP$, and observe something going wrong:
$$
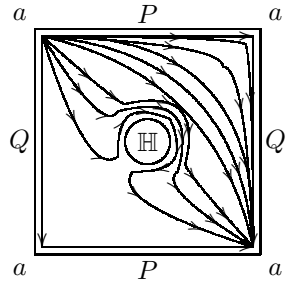
\xy
(15,15)*{ \mathbb{H}};
(15,15)*\xycircle(3,3){-};
{\ar (1,29)*{}; (29,29)*{}};
{\ar (29,29)*{}; (29,1)*{}};
(-2,-2)*{a};
(32,-2)*{a};
(32,32)*{a};
(15,32)*{P};
(32,15)*{Q};
(1,29)*{}; (28,28)*{} **\crv{(25,29)}?(.5)*\dir{>};
(28,28)*{};(29,1)*{} **\crv{(29,25)}?(.5)*\dir{>};
(1,29)*{}; (25,25)*{} **\crv{(25,28)}?(.5)*\dir{>};
(25,25)*{};(29,1)*{} **\crv{(28,25)}?(.5)*\dir{>};
(1,29)*{}; (29,1)*{} **\crv{(30,30)}?(.4)*\dir{>}?(.8)*\dir{>};
(1,29)*{}; (29,1)*{} **\crv{(27,27)}?(.4)*\dir{>}?(.8)*\dir{>};
(1,29)*{}; (29,1)*{} **\crv{(24,24)}?(.4)*\dir{>}?(.8)*\dir{>};
(1,29)*{}; (12,20)*{} **\crv{(11,20)}?(.4)*\dir{>}?(.8)*\dir{>};
(12,20)*{}; (20,12)*{} **\crv{(22.5,22.5)}?(.4)*\dir{>}?(.8)*\dir{>};
(20,12)*{}; (29,1)*{} **\crv{(20,11)}?(.4)*\dir{>}?(.8)*\dir{>};
(1,29)*{}; (11,19)*{} **\crv{(9,16)}?(.4)*\dir{>}?(.8)*\dir{>};
(11,19)*{}; (19,11)*{} **\crv{(22,22)}?(.4)*\dir{>}?(.8)*\dir{>};
(19,11)*{}; (29,1)*{}  **\crv{(16,9) }?(.4)*\dir{>}?(.8)*\dir{>};
(1,29)*{}; (11,13)*{}  **\crv{(7,10) }?(.4)*\dir{>}?(.8)*\dir{>};
(11,13)*{}; (18,18)*{}  **\crv{(11,19)&(13,20) }?(.4)*\dir{>}?(.8)*\dir{>};
(18,18)*{}; (13,11)*{}  **\crv{(20,13)&(19,11) }?(.4)*\dir{>}?(.8)*\dir{>};
(13,11)*{}; (29,1)*{}  **\crv{(10,7) }?(.4)*\dir{>}?(.8)*\dir{>};
(-2,15)*{Q};
{\ar (1,29)*{}; (1,1)*{}};
(15,-2)*{P};
{\ar (1,1)*{}; (29,1)*{}};
(-2, 32)*{a};
(0,0)*{}; (30,0)*{} **\dir{-};
(30,0)*{}; (30,30)*{} **\dir{-};
(30,30)*{}; (0,30)*{} **\dir{-};
(0,30)*{}; (0,0)*{} **\dir{-}
\endxy
$$
It's obvious that as we try to move the path $PQ$ continuously through intermediate paths
towards $QP$, we get stuck on the hole.
\smallskip

The fact that the order of the paths truly matters is what's meant by the adjective
`non-abelian,' which is the case for the fundamental group of $\cX$ (whatever that is, again).
It has quite surprising consequences for the counting problem.
As for the space $\cE$, for $\cX$ as well, it is a fact that any path that begins and ends
at $a$ is homotopic to a sequence of $P$'s and $Q$'s. However, in this case, we are {\em not} allowed
to change the order. For example, it turns out that
$$PPQP$$
is not homotopic to
$$PQPP.$$
So  let us count the paths of length 3.
We start again with three $P$'s:
$$PPP$$
then we change the last segment:
$$PPQ.$$
From here, we can proceed `lexicographically,' that is, list the paths according to the order they might
appear in a dictionary. We get afterwards
$$PQP;$$
$$PQQ;$$
$$QPP;$$
$$QPQ;$$
$$QQP;$$
$$QQQ;$$
making up a total of 8 paths of length 3. One sees already that the number is greater than what we obtained for
$\cE$. But it is probably not yet apparent how dramatic the difference is.
If you have a path of length $n$, it will look like
$$\Box \Box \Box \cdots \Box \Box$$
with either a $P$ or $Q$ in  the $n$ boxes.
But then there are two possibilities for every box, each of which
can be chosen independently of the other. From this perspective, another way of thinking about the count in the length three
case is that the first path could be $P$ or $Q$:
$$P \Box \Box;$$
$$Q \Box \Box.$$
The number of possibilities for the remaining two boxes will simply be the number of paths of length 2.
That is
$$\mbox{Number of paths of length 3}=2 \times \mbox{Number of paths of length 2}.$$
But then we can keep going iaccording to the same logic and say
$$\mbox{Number of paths of length 2}=2\times \mbox{Number of paths of length 1}$$
$$=2\times 2.$$
So
$$\mbox{Number of paths of length 3}=2\times 2\times 2=2^3=8.$$
This method of counting generalizes, so that
the number of paths of length $n$ is 2 multiplied to itself $n$ times, or
$$2^n.$$
Recall that in $\cE$, there were 11 paths of length $10$.
However, in $\cX$, the number is
$$2^{10}=1024.$$
When the length is increased to 100, the count in $\cE$ is
$$100+1=101,$$
while in $\cX$ it is
$$2^{100}.$$
It may not be obvious to you how staggeringly large this number is. If written out in the usual way,
it would be just around 1 followed by 30 zeros. I suspect you don't know any special name for
this number. But it might interest you to know that it is probably larger than the
number of atoms in the universe.
\smallskip

When we study equations, there are several different kinds of spaces that we associate to it, all of them useful in different ways. But the most naive one is called the {\em complex analytic space} of
the equation. This is, roughly put, the set of all complex number solutions of the equation, except
we have to change it a little bit to get a nice space. So we are trying to study the difficult
rational solutions by putting them into a very large set of solutions, large enough to have
a very intuitive geometric structure.

If we start with an equation like
$$x^3+y^3=1729,$$
then the analytic space has the shape of $\cE$.
For the equation
$$y^2=x^5 - 14x^4 + 65x^3 - 112x^2 + 60x$$
the associated space looks like two copies of $\cX$, glued along the edge of the holes that we've made.
That is, you would take two inner tubes, cut a hole in each of them, and then glue them together
along the hole, to create something like a Siamese twin doughnut.
\smallskip

What about the circle equation
$$x^2+y^2=1?$$
Its analytic space is a sphere, and there, {\em all} paths with a fixed endpoint $a$ are homotopic to each other, and
to the {\em null path} that just sits at $a$ and goes nowhere.
Therefore, the order of traversal will definitely not matter, accounting for the
infinitude of rational solutions.
\smallskip

I wish it were possible to say with some precision what these paths have to do with solutions.
If we're given an equation $f(x,y)=0 $, denote by $\cX_f$ its analytic space.
For any two points $a$ and $b$ on $\cX_f$, we have to consider the set
$$\P(a,b)$$ of all paths from $a$ to $b$, and how this set itself varies as we move the endpoints $a$ and $b$.
Wait, you say,  we weren't allowed to move the endpoints! That's correct, as far as homotopies are concerned.
But we are not moving the individual paths along homotopies. Rather, we are trying to keep track of how the
collection of all paths varies as we consider different possibilities for the endpoints.
\smallskip

In the space $\cX$, for example, we've indicated two different points $b$ and $c$,
two representative paths $P_1$ and $P_2$ belonging to $\P(a,b)$, and two paths
$Q_1, Q_2$ belonging to $\P(a,c)$. (Do  you see the trajectory of the path $P_2$?)
$$
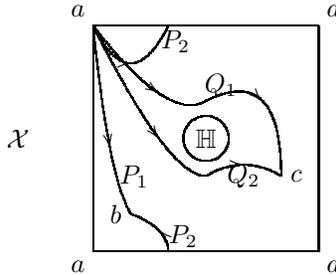
\xy
(-10,15)*{\cX};
(20,10)*{Q_2};
(15,10)*{}; (25,10)*{}  **\crv{(20,13) }?(.5)*\dir{>};
(0,30)*{}; (15,10)*{}  **\crv{(10,10) }?(.5)*\dir{>};
(17,22)*{Q_1};
(0,30)*{}; (15,20)*{}  **\crv{(10,17) }?(.5)*\dir{>};
(15,20)*{}; (25,10)*{}  **\crv{(25,25) }?(.5)*\dir{>};
(5.5,10)*{P_1};
(11,28)*{P_2};
(12,2)*{P_2};
(0,30)*{}; (5,5)*{}  **\crv{(2,10) }?(.5)*\dir{>};
(0,30)*{}; (10,30)*{}  **\crv{(5,20) }?(.5)*\dir{>};
(10,0)*{}; (5,5)*{}  **\crv{(10,3) }?(.5)*\dir{>}?;
(3,5)*{b};
(27,10)*{c};
(-2,-2)*{a};
(32,-2)*{a};
(32,32)*{a};
(15,15)*{ \mathbb{H}};
(15,15)*\xycircle(3,3){-};
(-2, 32)*{a};
(0,0)*{}; (30,0)*{} **\dir{-};
(30,0)*{}; (30,30)*{} **\dir{-};
(30,30)*{}; (0,30)*{} **\dir{-};
(0,30)*{}; (0,0)*{} **\dir{-}
\endxy
$$
To repeat, all the paths in $\P(a,b)$ have fixed endpoints $a$ and $b$, and homotopies between them
are not allowed to move the endpoints. Similarly, the paths in $\P(a,c)$ have fixed endpoints
$a,c$ and homotopies must leave those fixed. However, we are trying to understand how the  {\em collection}
$\P(a,b)$ differs from the collection $\P(a,c)$, or any $\P(a,x)$ as we now move $x$ around. This is somewhat mind-boggling, even for professional mathematicians, but we eventually manage to do it.
\smallskip

A  rational solution to the equation determines a very special kind of point on $\cX_f$, so that
when $a$ and $x$ are both rational solutions,   $\P(a,x)$ acquires a very large number of hidden symmetries associated to {\em Galois theory}, regardless of
the kind of space that $\cX_f$ is. The source of these symmetries is another very natural and highly symmetric space
$\cX_f^{et}$, unfortunately  hard to visualize. Our paths also move inside $\cX_f^{et}$ in very structured ways, even though
the space itself has something of a `fractal' nature,
and intersects the intuitive space $\cX_f$ in a sparse collection of points, the so-called
`algebraic points.'
When  symmetries are taken into
account,  $\P(a,x)$ is actually capable of distinguishing the points.
That is, $\P(a,b)$ and $\P(a,c)$ will look quite different
whenever $b$ and $c$ are different. In the situation where $\cX_f$ has non-abelian fundamental group,
the consequent complexity of any one of the $\P(a,x)$ creates a severe tension. The conflict between
the high degree of symmetry and this complexity can only be resolved by an extreme rarity of such
structures. That is, only finitely many $\P(a,x)$'s should be possible.
 I should admit, however, that even with a powerful result like Faltings' theorem in hand,  it's not yet clear
how uniform an explanation this will turn out to be. Of course, we are right now trying
to clarify this issue as part of a detailed programme to make {\em constructive use} of this complexity.

\end{flushleft}

\end{document}